# PROPERTIES OF SUMS OF SOME ELEMENTARY FUNCTIONS AND MODELING OF TRANSITIONAL AND OTHER PROCESSES

Yuri K. Shestopaloff

*Abstract*: The article presents mathematical generalization of results which originated as solutions of practical problems, in particular, the modeling of transitional processes in electrical circuits and problems of resource allocation. However, the presented findings have broader meaning and can be used for approximation of transitional and other processes in different areas of science and technology. We present discovered properties of sums of polynomial, power, and exponential functions of one variable. It is shown that for an equation composed of one type of function there is a corresponding equation composed of functions of the other type. The number of real solutions of such equations and the number of characteristic points of certain appropriate corresponding functions are closely related. In particular, we introduce a method similar to Descartes' Rule of Signs that allows finding the maximum number of real solutions for the power equation and equation composed of sums of exponential functions. The discovered properties of these functions allow us to improve the adequacy of mathematical models of real phenomena.

## 1. Formulation of the problem

Probably, everybody is familiar with residual oscillation effects in mechanical systems. For instance, when we start or turn off car engine, it first moves (oscillates) before it comes to rest or steady state. When one stops a centrifuge juicier, the juicier can make few violent movements before it stops. Similar oscillation effects, observed during the transition from one functional state to another, experience many other mechanical systems and mechanisms.

All electrical systems and devices include in one form or another transitional electrical processes. Modern electronic devices

In general, it is easier to model steady states of dynamic systems, such as electrical or mechanical, while it is more difficult to find an adequate description for transitional processes. On the other hand, knowledge of characteristics of transitional processes is important, at least for several reasons. Transitional states often impose stressful conditions to mechanisms and systems, both in terms of acting forces and amplitudes. Another reason is the need to accurately estimate when transition phase ends, based on some criteria, so that we can consider the system's state as stationary. Examples of such electrical systems are presented in [1].

Often, transitional processes are described by attenuating or increasing exponential function, such as, for instance, single *RC* series circuit, or by sums of exponential functions. In case of electrical systems composed of many elements, the general solution for transition process may be difficult to find. In this case, the required characteristics of transitional process, for instance, time and number of oscillations before the system becomes stable, can be accurately found based on approximation of transitional process as sums of exponential or other functions. In particular, this approach has been used to find the transition time and maximum signal spike in electrical circuits that process short impulses. Then, we used mathematical results with regard to properties of sums of elementary functions, presented in this paper, such as the theorem about the maximum number of oscillations of a sum of exponential signals. The discovered mathematical properties much facilitated the solution of the problem, which was the reduction of transitional time and substantial diminishing of value of signal spike by virtue of optimal spatial configuration of components, and characteristics of electrical impulses.

Also, we discovered that the results can be beneficially applied to study of similar problems in other areas, in particular, propagation of arbitrary polarized waves in electrically anisotropic and highly absorbing mediums, when the wideband signal can be approximated by discrete attenuating signals with distinct polarizations.

The discovered properties of sums of elementary functions also have found application in the area of optimal allocation of resources. In particular, the IRR equation, which is the foundation of many mathematical methods used for analysis and allocation of financial resources, is a particular case of power equation which properties we consider in this paper. Application of presented theorems allowed to strictly define the business meaning of solutions of IRR equation [2], and offer several computational approaches how to find the number of solutions for this equation in a general case, and what solution has to be chosen for practical applications [3]. In some practically important scenarios, it is possible to unambiguously define the domain, in which the required solution is located [3].

In this article, the findings that were previously discovered as solutions of practical problems mentioned above are presented in a generalized form.

## 2. Introducing a corresponding equation for a polynomial

Mathematical modeling of real phenomena is often based on usage of elementary functions and their sums [1, 4-7]. In this paper, we show that the properties of sums of some elementary functions are related.

For the polynomial equation, we introduce a corresponding equation, composed of the following sum of exponential functions. Let us consider the polynomial.

$$\sum_{j=0}^{N} C_j y^{T_j} = 0 \tag{1}$$

where $C_j$ are real constants. We assume that $T_j$ are integers, because this is a polynomial equation; $y$ is a variable; zero value of $T_j$ is permitted in order to include a free term into (1).

The maximum possible number of real solutions is defined by the number of sign changes $N$ of consecutive terms in (1) written in the descending order of powers, according to Descartes' Rule of Signs. In our case, the maximum number of real roots is $(2N+1)$, which may include multiplicities (zero solution is counted once).

We assume that $y > 0$, so that we can do the Taylor series expansion of the power functions [5] at the point $T_j = 0$ as follows. (In our case, these are polynomials, because $T_j$ are integers. Power functions are similar to polynomials but they have the real powers.)

$$\sum_{j=0}^{N} C_j \sum_{k=0}^{\infty} \left( \frac{T_j^k \ln^k(y)}{k!} \right) = \sum_{k=0}^{\infty} \frac{\ln^k(y)}{k!} \sum_{j=0}^{N} C_j T_j^k \tag{2}$$

Then, we can rewrite (1) in the following way.

$$\sum_{k=0}^{\infty} \frac{\ln^k(y)}{k!} \sum_{j=0}^{N} C_j T_j^k = 0 \tag{3}$$

The logarithmic function is a one to one function. So, (2) is an equivalent transformation for $y > 0$. Substituting $x = \ln(y)$, we can rewrite (3) as follows.

$$\sum_{k=0}^{\infty} x^k \sum_{j=0}^{N} \frac{C_j T_j^k}{k!} = \sum_{k=0}^{\infty} a_k x^k = 0 \tag{4}$$

where $a_k = \sum_{j=0}^{N} \frac{C_j T_j^k}{k!}$.

So, instead of the original equation (2), for $y > 0$, we can analyze the solutions of equation (4), which can be rewritten as follows.

$$\sum_{k=0}^{\infty} a_k x^k = 0 \tag{5}$$

For $y < 0$, we can do similar transformations, rewriting (1) as follows.

$$-\sum_{j=0}^{N} C_j y^{T_j} = \sum_{j=0}^{N} (-1)^{T_j+1} C_j (-y)^{T_j} = 0 \tag{6}$$

Consequently, we can rewrite (6), similar to (2), in the following way.

$$\sum_{k=0}^{\infty} \frac{\ln^k(-y)}{k!} \sum_{j=0}^{N} (-1)^{T_j-1} C_j T_j^k = 0 \tag{7}$$

Instead of the original equation, we can analyze the number of solutions of the following equation.

$$\sum_{k=0}^{\infty} a_k^{(-)} x^k = 0 \tag{8}$$

where

$$a_k^{(-)} = \sum_{j=0}^{N} (-1)^{T_j-1} C_j T_j^k \tag{9}$$

Note that despite the increased range of argument in (3), (logarithm function can be both positive and negative), the maximum number of solutions remains the same, because the logarithmic function does a one to one transformation. So, when we combine (5) and (8), we will get the same maximum number of sign changes as for the polynomial equation, that is $(2N+1)$.

## 3. Relationship between the number of sign changes of exponential sums and polynomial's solutions

Note that transformations from (1) to (5), and from (6) to (8) are mathematically *equivalent*. So, the number of solutions of the polynomial equation (1) is equal to the number of solutions of equation (5), whose coefficients are composed of sums of exponential functions. The maximum number of real solutions in (1) is defined by the number of sign changes counted according to the Descartes' Rule of Signs. For positive *y*, this means that the maximum possible number of real solutions is *N*.

Let us consider (5), that is $\sum_{k=0}^{\infty} a_k x^k = 0$, and coefficients $a_k = \sum_{j=0}^{N} \frac{C_j T_j^k}{k!}$. Suppose that $T_p = \max\{T_j\}, j = \{0,..,N\}$. Then, we can do the following transformation.

$$a_k = \sum_{j=0}^{N} \frac{C_j T_j^k}{k!} = \frac{1}{k!} C_p T_p^k \left( 1 + \sum_{j=0, j\neq p}^{N} \frac{C_j}{C_p} \times \left( \frac{T_j}{T_p} \right)^k \right) \tag{10}$$

It follows from (10) that

$$\begin{aligned} \lim_{k\to\infty} a_k &= \lim_{k\to\infty} \frac{1}{k!} C_p T_p^k \left( 1 + \sum_{j=0, j\neq p}^{j=N} \frac{C_j}{C_p} \times \left( \frac{T_j}{T_p} \right)^k \right) = \\ &= \lim_{k\to\infty} \frac{1}{k!} C_p T_p^k = 0 \end{aligned} \tag{11}$$

However, what is important to us is the change of sign of $a_k$. The term $C_p T_p^k$ preserves its sign for all *k*. The term $\left( 1 + \sum_{j=0, j\neq p}^{j=N} \frac{C_j}{C_p} \times \left( \frac{T_j}{T_p} \right)^k \right)$ is such that

$$\lim_{k\to\infty} \left( 1 + \sum_{j=0, j\neq p}^{j=N} \frac{C_j}{C_p} \times \left( \frac{T_j}{T_p} \right)^k \right) = 1$$

It means that $\sum_{j=0, j\neq p}^{j=N} \frac{C_j}{C_p} \times \left( \frac{T_j}{T_p} \right)^k = o(1)$, and for any small positive value $\varepsilon$ such a number *K* always exists that for all $k > K$ the following inequality holds true

$$\left| \sum_{j=0, j \neq p}^{j=N} \frac{C_j}{C_p} \times \left( \frac{T_j}{T_p} \right)^k \right| < \varepsilon$$

In fact, in order to preserve the sign of the whole expression

$$C_p T_p^k \left( 1 + \sum_{j=0, j \neq p}^{j=N} \frac{C_j}{C_p} \times \left( \frac{T_j}{T_p} \right)^k \right)$$

in (11), we need even weaker condition, namely $\left| \sum_{j=0, j \neq p}^{j=N} \frac{C_j}{C_p} \times \left( \frac{T_j}{T_p} \right)^k \right| < 1$. So, we proved that such a number *K* always exists that for all $k > K$ the sign of all $a_k$ remains the same. So, *all* sign changes of coefficients $a_k$ are always confined within the *finite* interval $\{0, K\}$.

We assume *K* finite but large enough, so that the number of changes of algebraic sign of coefficients $a_k$ is equal to *N,* because (5) is an equivalent transformation of (1). (We remind that *N* is the number of sign changes of the original polynomial equation, which, in case of $y > 0$, has a maximum of *N* solutions.) Recall that (5) converges to (1) as the Taylor series representation of (1). Our large number *K* affects only the term from which we begin to analyze the sign changes, but not the Taylor series representation (5), in which the number of terms remains infinite.

Given the definition of $a_k$, it means that the number of sign changes of the series composed of the sum of exponential functions $S_k = \sum_{j=0}^{N} \frac{C_j T_j^k}{k!}$, where $k = \{0, \infty\}$, cannot be more than $N$ too. It gives the rise to an assumption that probably the following equation, whose left part represents the sum of exponential functions, has the same number of solutions with regard to unknown value of *k* (if we assume that *k* is real, then the function is continuous, which means that each change of sign is accompanied by the abscissa intersection).

$$\sum_{j=0}^{N} C_j T_j^k = 0 \tag{12}$$

However, the objection can be that if *k* is continuous, then (12) can have more solutions than the number of sign changes in series $S_k$. We have yet to prove that this is not the case and our result is valid for real value of *k* in (12). However, before doing this, we have to generalize our finding for the real powers $T_j$, which will be the topic of next section.

Thus, we have an equivalency of the maximum number of solutions of the original polynomial equation (1) and the number of sign changes of the series composed of sums of exponential functions $S_k = \sum_{j=0}^{N} \frac{C_j T_j^k}{k!}$. Note that the presence of the factorial does not influence the number of sign changes, which will be the same for the sum of exponential functions $S = \sum_{j=0}^{N} C_j T_j^k$ .

We should make the following note about the free term in (1) and (5). The sums may also have a free term depending on the value of parameters. Geometrically, an addition of the free term means the displacement of the graph of function $F(y) = \sum_{j=0}^{N} C_j y^{T_j}$ along the ordinate axis, which in some instances can bring additional (real) solution. This note complies with the Descartes' Rule of Signs, because the addition of the free term can add extra change of sign. So, we assume that the free term can present both in the polynomial equation and its corresponding equation composed of the sum of appropriate exponential functions. We can formulate the obtained results as the following Lemma.

**Lemma 1:**

Each polynomial equation $\sum_{j=0}^{N} C_j y^{T_j} = 0$ has a corresponding series composed of the sum of exponential functions, defined as follows: $S_k = \sum_{j=0}^{N} C_j T_j^k$ , $k = 0,1,..\infty$. The reverse is also true, that is each series of sums of exponential functions with positive integer bases has an associated unique polynomial equation. The maximum number of solutions of the polynomial equation is equal to the maximum number of algebraic sign changes of its corresponding series composed of sums of exponential functions.

## 4. Generalization of Lemma 1 for real powers

Let us consider the power equation with rational powers. In this case, (1) can be rewritten as follows.

$$\sum_{j=0}^{N} C_j y^{\frac{m_j}{n_j}} = 0 \qquad (13)$$

Substituting $z = y^{\frac{1}{n_j}}$, we can transform (13) into an equation equivalent to polynomial (1).

$$\sum_{j=0}^{N} C_j z^{m_j} = 0 \qquad (14)$$

According to Lemma 1, this equation has the same maximum number of solutions as the number of sign changes of series. So, in case of rational powers, the number of real solutions is defined by the number of sign changes of series composed of sums of exponential functions, which in turn is equal to the number of sign changes of consecutive terms of the power equation written in the descending order of powers. Note that the substitution $z = y^{\frac{1}{n_j}}$ can bring an additional negative solution. However, when we consider the total number of solutions for the negative and positive argument separately, then the total number of solutions will remain the same.

Let us consider the irrational number $T_j$. In this case, we can always approximate it by the rational number $\frac{m_{pj}}{n_{pj}}$ such that

$$\left| T_j - \frac{m_{pj}}{n_{pj}} \right| < \delta \qquad (15)$$

where $\delta > 0$, so that

$$\lim_{p \to \infty} \frac{m_{pj}}{n_{pj}} = T_j$$

The equation

$$\sum_{j=0}^{N} C_j y^{\frac{m_{pj}}{n_{pj}}} = 0$$

has a number of solutions defined by the number of sign changes of consecutive terms, which means that its corresponding series composed of the appropriate sums of exponential functions also has the same number of sign changes.

Given (15) and the continuity and smoothness of the considered functions (its first derivative is also continuous function), it means that for any arbitrary large domain of $y$, and for any small value $\varepsilon > 0$ such number $\frac{m_{pj}}{n_{pj}}$ and $\delta$ defined in (15) exist that

$$\left| \sum_{j=0}^{N} C_j y^{T_j} - \sum_{j=0}^{N} C_j y^{\frac{m_{pj}}{n_{pj}}} \right| < \varepsilon$$

We can find value $\delta$ as follows. Let us do the following transformations in the above inequality.

$$\left| \sum_{j=0}^{N} C_j y^{T_j} \left( y^{T_j - \frac{m_{pj}}{n_{pj}}} - 1 \right) \right| \le \sum_{j=0}^{N} \left| C_j y^{T_j} \left( y^{T_j - \frac{m_{pj}}{n_{pj}}} - 1 \right) \right| \le \sum_{j=0}^{N} \left| C_j y^{T_j} \left( y^{\delta} - 1 \right) \right| =$$

$$\left| y^{\delta} - 1 \right| \sum_{j=0}^{N} \left| C_j y^{T_j} \right| \le K \left| 1 + \frac{\delta \ln y}{1!} - 1 \right| \sum_{j=0}^{N} \left| C_j y^{T_j} \right| \le K\delta \left| \ln y \right| \sum_{j=0}^{N} \left| C_j y^{T_j} \right| < \varepsilon \qquad (16)$$

Here, we use the Taylor series expansion of function $y^{\delta}$ and its convergence, which means that such value $K > 1$ always exists that

$$y^{\delta} < K \left| 1 + \frac{\delta \ln y}{1!} \right|$$

Solving (16) with regard to $\delta$, we will find

$$\delta < \frac{\varepsilon}{K \left| \ln y \right| \sum_{j=0}^{N} \left| C_j y^{T_j} \right|}$$

So, we proved that for any small positive value $\varepsilon$ such a value $\delta$ always exists that the power function with the rational powers approximate the power function with irrational exponents with any required accuracy. Note that $\delta$ depends on $y$, which can be infinitely large. However, we are interested in the $\varepsilon$ accuracy approximation in the *finite* domain of $y$, which contains all solutions of the power equation. These solutions are always contained within the finite domain, because the sum of power functions always goes asymptotically to minus or plus infinity when $y \to \pm\infty$ [5]. So, we can always approximate the sum of power functions with irrational powers by the sum of appropriate power functions with rational powers with any required accuracy in the finite domain that includes all solutions.

Consequently, the power equation $\sum_{j=0}^{N} C_j y^{T_j} = 0$ with irrational powers will have the same maximum number of solutions as the limit of the maximum number of solutions of the appropriate equation $\sum_{j=0}^{N} C_j y^{\frac{m_j}{n_j}} = 0$ with rational powers, when $\delta \to 0$. Graphically, it means that the graphs of functions $F_T = \sum_{j=0}^{N} C_j y^{T_j}$ and $F_R = \sum_{j=0}^{N} C_j y^{\frac{m_j}{n_j}}$ coincide at the limit, when $\delta \to 0$, see Fig. 1. In other words, $\lim_{\delta \to 0} F_R = F_T$ for any arbitrary large domain of $y$ that includes all solutions of equation $\sum_{j=0}^{N} C_j y^{T_j} = 0$.

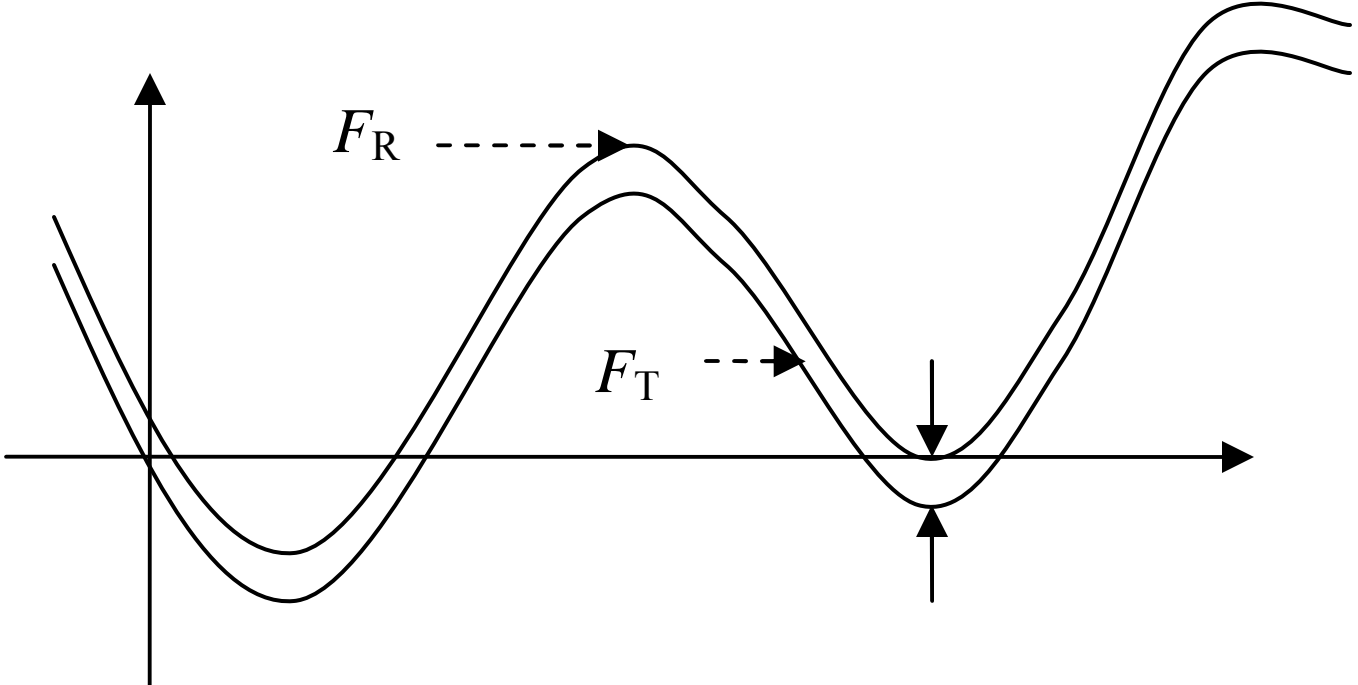

Fig. 1. Approximating the power function with irrational powers by the power function with rational powers.

So, we have a correspondence between the maximum number of solutions of the power equation with rational powers (the maximum number of solutions remains the same regardless of the value of powers, because the number of terms and consequently possible sign changes remains the same) and the maximum number of solutions of the power equation with irrational powers. So, the power equation with irrational powers has the same maximum number of solutions as the power equations with rational powers that approximate it. Thus, in case of a power equation with irrational powers, the maximum number of solutions is also defined by the number of sign changes of its consecutive terms written in the descending order of powers.

Similarly, we can show that each power equation with irrational powers has a corresponding series of sums of exponential functions, so that the number of sign changes of this series $S = \sum_{j=0}^{N} C_j T_j^k$ is equal to the maximum number of real solutions of a power equation. The proof repeats all steps, when we did this for the power equation with irrational powers. First, we create the series of exponential functions corresponding to the power function with irrational powers, that is

$$S_{Tk} = \sum_{j=0}^{N} C_j T_j^k$$

Then, we do the same for the power function with rational powers that approximates the original power function.

$$S_{Rk} = \sum_{j=0}^{N} C_j \left(T_{Rj}\right)^k$$

Here, $T_{Rj}$ denotes the approximating rational numbers.

We proved already that the power functions with irrational and rational powers can be close with any required accuracy by choosing the appropriate value of the approximating powers. These two sets of irrational and rational powers constitute accordingly two sets of series of sums of exponential functions introduced above. The difference between them can be found as follows.

$$\left|S_{Tk} - S_{Rk}\right| = \sum_{j=0}^{N} \left|C_j\right| T_{Rj}^k \left| \left(\frac{T_j}{T_{Rj}}\right)^k - 1 \right|$$

We know that for any small value of $\delta$ such a rational number $T_{Rj}$ exists that $\left|T_{Rj} - T_j\right| < \delta$.

So, we can rewrite the above expression as follows.

$$\left|S_{Tk} - S_{Rk}\right| \le \sum_{j=0}^{N} \left|C_j T_{Rj}^k\right| \left| \left(\frac{T_{Rj} + \delta}{T_{Rj}}\right)^k - 1 \right| = \sum_{j=0}^{N} \left|C_j T_{Rj}^k\right| \left| (1 + \frac{\delta}{T_{Rj}})^k - 1 \right| =$$

$$= \sum_{j=0}^{N} \left|C_j T_{Rj}^k\right| \left| 1 + \frac{k\delta}{T_{Rj}} + \ldots - 1 \right|$$

We proved earlier that all sign changes occur within the finite range of *k*. So, we can assume that *k* is finite, $k < K$. Then, we can find the limit

$$\lim_{\delta \to 0} \left|S_{Tk} - S_{Rk}\right| = \sum_{j=0}^{N} \left|C_j T_{Rj}^k\right| \lim_{\delta \to 0} \left| \frac{K\delta}{T_{Rj}} + \frac{K(K-1)\delta^2}{2!T_{Rj}^2} + \ldots \right| = 0$$

for any $k < K$. In other words, for any small value $\varepsilon$ such value $\delta$ exists that for all $\left|T_{Rj} - T_j\right| < \delta$ we have $\left|S_{Tk} - S_{Rk}\right| < \varepsilon$ for any $k < K$, which means that two series of sums of exponential functions can be infinitely close. (Note that we found the total limit as the sum of limits for each term of the sum; each such limit is equal to zero. It means that each term of the sum of exponential functions with rational powers converges to the appropriate term in the sum of exponential functions with rational powers.) The convergence of two series of sums of exponential functions and corresponding terms in these sums means that these series have equal number of sign changes.

So, we conclude that the number of solutions of equation $\sum_{j=0}^{N} C_j y^{T_j} = 0$ with real $T_j$ , and the number of sign changes in its corresponding series of sums of exponential functions are the same. Thus, we can generalize Lemma 1 as follows.

**Theorem 1***:*

*(correspondence of power functions and series of sums of exponential functions)*

Each power equation $\sum_{j=0}^{N} C_j y^{T_j} = 0$ with real powers $T_j$ has a corresponding unique series composed of the sum of exponential functions, defined as follows: $S_k = \sum_{j=0}^{N} C_j T_j^k$ . The reverse is also true, that is each series of sums of exponential functions with positive real bases has a corresponding unique power function and equation. The maximum number of real solutions of the power equation, in which all terms are real, is defined by the number of sign changes of consecutive terms written in the descending order of power; which is also equal to the number of sign changes of the corresponding series composed of sums of exponential functions.

## 5. Equivalency of solutions of the power equation and its corresponding equation composed of the sum of exponential functions

In this section, we will prove that the maximum number of real solutions of power equation is equal to the number of solutions of an equation composed of sum of appropriate exponential functions. In other words, that the equation $\sum_{j=0}^{N} C_j T_j^k = 0$ has the same number of real solutions as its corresponding power equation $\sum_{j=0}^{N} C_j y^{T_j} = 0$, and vice versa.

Let us consider the following sum of exponential functions and its transformation.

$$\sum_{j=0}^{N} C_j T_j^{\frac{k}{n_j}} = \sum_{j=0}^{N} C_j \left( T_j^{\frac{1}{n_j}} \right)^k = \sum_{j=0}^{N} C_j \left( T_{aj} \right)^k \tag{19}$$

where $T_{aj} = T_j^{\frac{1}{n_j}}$ .

Based on theorem 1, there is a unique power equation corresponding to the sum of exponential functions on the right side of (19), which is defined as follows.

$$\sum_{j=0}^{N} C_j y^{T_{aj}} = 0 \qquad (20)$$

According to the theorem 1, this equation has a maximum of $2N+1$ real solutions with multiplicities ($N$ is the number of sign changes of consecutive coefficients, for the positive argument and accordingly the number of solutions for $y \geq 0$); and the series $S_k = \sum_{j=0}^{N} C_j \left(T_{aj}\right)^k$ has $N$ maximum number of sign changes. (Creating a similar sum of exponential functions for negative $y$ based on (8) will compliment the total number of sign changes to $(2N+1)$.) Taking into account (19), this also means that the series

$$S_k = \sum_{j=0}^{N} C_j T_j^{\frac{k}{n_j}} \qquad (21)$$

has the same $N$ number of sign changes. Note that we did not impose any limitation to the value of $n_j$, except that it is positive. So, the power in (21) can be any *rational* number, while the maximum number of sign changes remains the same regardless of the value $n_j$. It is known from calculus that any irrational number $R$ can be approximated by a rational number with any required accuracy $\delta_R$, so that

$$\left| R - \frac{m}{n} \right| < \delta_R \qquad (22)$$

Given the continuity of exponential functions, this means that (similar to what we did in (14) ) for any small number $\varepsilon$, such rational numbers always exist that

$$\left| \sum_{j=0}^{N} C_j T_j^{R_j} - \sum_{j=0}^{N} C_j T_j^{\frac{m_j}{n_j}} \right| < \varepsilon \qquad (23)$$

We can transform (23) and find the limit, taking into account (22), as follows.

$$\lim_{R_j - \frac{m_j}{n_j} \to 0} \left| \sum_{j=0}^{N} C_j T_j^{\frac{m_j}{n_j}} \left(T_j^{R_j - \frac{m_j}{n_j}} - 1\right) \right| = 0 \qquad (24)$$

The result (24) means that the sum of exponential functions with irrational powers can be approximated by the sum of exponential functions with rational powers with any required accuracy. Hence, the maximum number of sign changes of the sum of exponential functions with irrational powers is the same as the number of sign changes for the approximating sum of exponential functions with rational powers. So, we proved that in the series of sums of exponential functions the maximum number of sign changes remains the same for any ascending series of arguments $\{x_k\}, 0 \le x_k < \infty$, where $x_k$ is any non-negative real number, not only integers, as we proved previously. So, we can generalize the series $S_k = \sum_{j=0}^{N} C_j T_j^k$ as follows.

$$S(x_k) = \sum_{j=0}^{N} C_j T_j^{x_k} \tag{25}$$

It also means that the graph of function

$$S(x) = \sum_{j=0}^{N} C_j T_j^{x} \tag{26}$$

intersects the abscissa maximum the same number of times as is the number of sign changes. In other words, this is the number of solutions that the following equation, which is composed of the sum of exponential functions, has.

$$\sum_{j=0}^{N} C_j T_j^{x} = 0 \tag{27}$$

The note should be made with regard to the existence of such a finite number *X* that for all $x > X$ the sign of $\sum_{j=0}^{N} C_j T_j^x$ remains the same. The result follows directly from the fact that when $x \to \infty$ the function with the largest base dominates, so that the rest of summands is o-small of this dominant function. Similar, when $x \to -\infty$, the summand with the smallest base dominates [5].

Thus, we came to the result that can be formulated as the following Theorem.

**Theorem 2**

*(correspondence of exponential and power functions and equations):*

The power equation $\sum_{j=0}^{N} C_j y^{T_j} = 0$ and its corresponding equation composed of the sum of exponential functions $\sum_{j=0}^{N} C_j T_j^x = 0$ have the same maximum number of real solutions defined by the number of sign changes, when the power equation is written in the descending order of powers. Alternatively, this is also the number of sign changes in the corresponding sum of exponential functions written in the descending order of bases of exponential functions.

## 6. Logarithmic functions

In this short section, we would like to pay attention to the fact that the sum of logarithmic functions is a monotonic function. This property explains the behavior of many natural phenomena, in particular the growth processes understood in a broad sense [6, 7]. A logarithmic function is the mathematical inverse of exponential function. If we write the exponential function in a form $x = Ct^{y}$, then $y = \frac{1}{\ln(tC)} \ln x$ for $x > 0$. Graphically, the logarithmic function can be obtained from the exponential function by a mirror reflection relative to axis $y = x$, see Fig. 2.

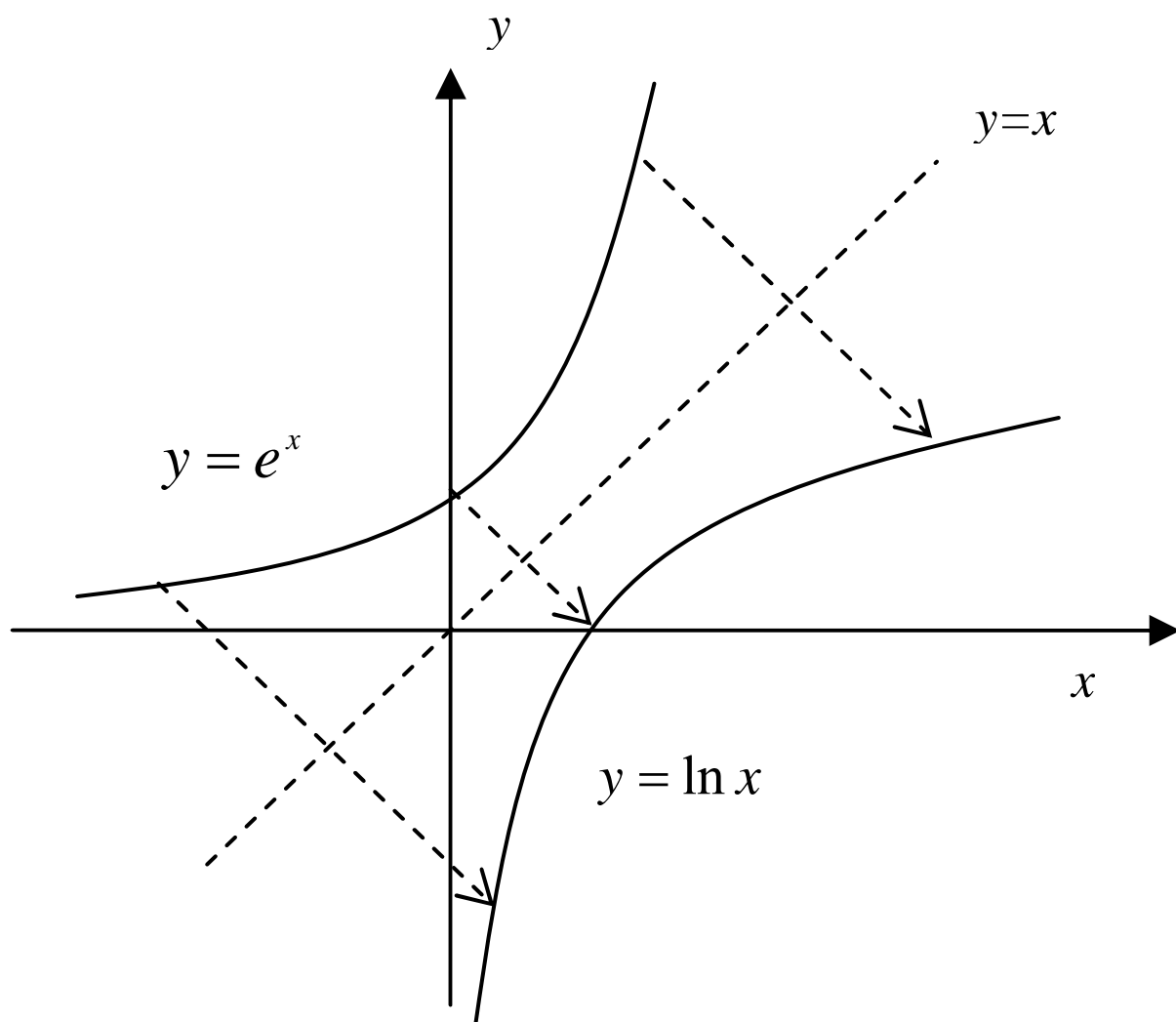

Fig. 2. Logarithmic function as an inverse of exponential function.

The sum of logarithmic functions $y(x) = \sum_{j} C_j \log_{a_j} x$ can be transformed into one logarithmic function with an arbitrary base, say to a natural logarithm. We can prove this as follows.

$$y(x) = \sum_{j} C_j \log_{a_j} x = \sum_{j} \frac{C_j}{\ln a_j} \ln x = \ln x \sum_{j} \frac{C_j}{\ln a_j} = C_0 \ln x \tag{28}$$

where $C_0 = \sum_j \frac{C_j}{\ln a_j}$ is a constant.

We emphasized the monotonicity of the sum of logarithmic functions with regard to their application in modeling and mathematical description of many natural phenomena, so that it should not be viewed as a purely abstract mathematical property, but as a direct relationship of logarithmic and exponential functions with Nature.

## 7. The theorem about the correspondence of sums of function and its application

The importance of the proved theorems is that previously these types of functions and equations have been considered much unrelated. Given the significant role of exponential equations in the modeling of natural phenomena, this theorem delivers a very powerful tool to handle problems in which these types of functions are involved.

From now on, the analysis of certain types of functions and equations is not restricted by the range of these functions, but many more analytical methods can be used that have been developed for other types of corresponding functions. For instance, exponential functions and equations have substantially less analytical methods than the polynomial equations and functions, so that these results can be transformed to analysis of exponential functions.